\documentclass{amsart}
\usepackage{geometry} 
\usepackage{url}
\usepackage{amsthm}
\usepackage{amssymb}
\usepackage{latexsym}
\usepackage{amsfonts}
\usepackage{amsmath}
\usepackage{amstext}
\usepackage{accents}
\usepackage{cancel}
\usepackage{eucal}
\usepackage{amscd}
\usepackage{hyperref}
\usepackage{graphicx}
\usepackage[shortalphabetic]{amsrefs}
\newtheorem{theorem}{Theorem}
\newtheorem{lemma}[theorem]{Lemma}
\newtheorem{corollary}[theorem]{Corollary}
\newtheorem{proposition}[theorem]{Proposition}
\theoremstyle{definition}
\newtheorem*{definition}{Definition}
\newtheorem*{conj}{Conjecture}
\newcommand{\R}{\mathbb{R}}

\newcommand{\RR}{\mathbb{R}}
\renewcommand{\SS}{\mathbb{S}}
\renewcommand{\S}{\mathbb{S}}

\newcommand{\SK}{\mathcal{K}}
\newcommand{\SF}{\mathcal{F}}
\newcommand{\SC}{\mathcal{C}}
\newcommand{\SD}{\mathcal{D}}
\newcommand{\SQ}{\mathcal{Q}}
\newcommand{\ip}[2]{\ensuremath{\left\langle #1, #2 \right\rangle}}

\newcommand{\vn}[1]{\ensuremath{\lVert#1\rVert}}

\begin{document}

\title[Solitons for the curve diffusion flow]{The shrinking figure eight and other solitons for the curve diffusion flow}
\author[M. Edwards]{Maureen Edwards}
\address{Institute for Mathematics and its Applications, University of Wollongong}
\email{maureen@uow.edu.au}
\author[A. Gehardt-Bourke]{Alexander Gerhardt-Bourke}
\address{Institute for Mathematics and its Applications, University of Wollongong}
\email{agb526@uowmail.edu.au}
\author[J. McCoy]{James McCoy*}\thanks{*: Corresponding author}
\address{Institute for Mathematics and its Applications, University of Wollongong}
\email{jamesm@uow.edu.au}
\author[G. Wheeler]{Glen Wheeler}
\address{Institute for Mathematics and its Applications, University of Wollongong}
\email{glenw@uow.edu.au}
\author[V.-M. Wheeler]{Valentina-Mira Wheeler}
\address{Institute for Mathematics and its Applications, University of Wollongong}
\email{vwheeler@uow.edu.au}
%\author[G. Williams]{Graham Williams}
%\address{Faculty of Law, Humanities and The Arts, University of Wollongong}
%\email{ghw@uow.edu.au}
\thanks{The research of the second author was supported by a summer vacation
scholarship of the Australian Mathematical Sciences Institute.  The research of
the third, fourth and fifth authors was supported by Discovery Project grant
DP120100097 of the Australian Research Council.}
\begin{abstract}
In this article we investigate the dynamics of special solutions to the surface diffusion flow of idealised ribbons.
This equation reduces to studying the curve diffusion flow for the profile curve of the ribbon.
We provide: (1) a complete classification of stationary solutions; (2)
qualitative results on shrinkers, translators, and rotators; and (3) an
explicit parametrisation of a shrinking figure eight curve.
\end{abstract}

\maketitle

\section{Introduction} \label{S:intro}

The curve diffusion flow is a particular evolution law governing the changing
shape of a curve in time.  It is the one-intrinsic-dimension analogue of the
surface diffusion flow for evolving surfaces.  These flows are gradient flows
of a certain energy functional in a certain function space, but not the more
common $L^2$-gradient flow such as the curve shortening flow (corresponding to
decreasing the length functional for a curve), the mean curvature flow
(decreasing the area functional for a surface), or the Willmore flow
(decreasing the elastic energy for a curve or the Willmore energy for a
surface). 

Surface diffusion flow was first described by W. W. Mullins in 1956 as a model
for the development of grooving at grain boundaries of heated crystal
structures \cite{Mullins}.  The flow is also related to the Cahn-Hilliard
equation, which describes phase separation and coarsening in the quenching
process of binary alloys \cite{CENC}. For further analysis and historical
remarks on the surface diffusion flow and its constrained variations we refer
the interested reader to \cites{MWW,W1,W2}.

The surface diffusion flow is a fourth-order nonlinear parabolic system of
equations that, at each time $t$, moves points $p$ on the solution surface $f$
perpendicularly to the surface with speed $(\Delta H)(p)$. Here $\Delta$ is the
Laplace-Beltrami operator and $H$ is the mean curvature, which is given by $H =
\ip{\Delta f}{\nu}$. The inner product $\ip{\cdot}{\cdot}$ is the standard
dot product in $\RR^3$.

Surface diffusion flow is the geometric analogue of the classical parabolic
clamped plate equation, which expresses the dynamics of a flat plate realised
as the graph of a function over a domain in $\RR^2$.
The curve diffusion flow is similarly the geometric analogue of the classical
parabolic Euler-Bernoulli model for a clamped rod, modelled as the graph of a
function over an interval in $\RR$.
When the domain is no longer flat, as is the case for a clamped arc segment of
a circle in the plane or cylindrical ring in space, the geometry of the problem
must enter into the model \cites{ESK,ES}.
The surface or curve diffusion flow can capture the notion of stability for an
elastic plate or rod under stress respectively.
Indeed, if the flow returns all small perturbations of a given initial
configuration to that initial configuration without exception, then this
configuration is called stable \cite{N}.

If the evolving surfaces are cylindrical, as would be the case for idealised
closed ribbons with their width extended to infinity, then this remains so
under the flow.
In this case the flow equation reduces to a one-dimensional parabolic equation
for the profile curve $\gamma$.
For the surface diffusion flow, this equation is exactly the curve diffusion
flow.
This reduction is also performed in \cite{ES}.

\begin{figure}
\includegraphics[width=3cm]{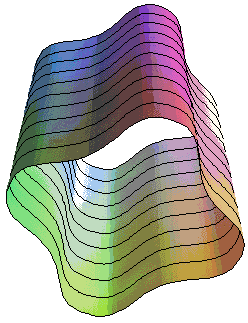}
\hspace{2cm}
\includegraphics[width=3cm]{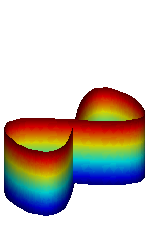}
\caption{Pieces of two cylindrical ribbons, immersed via $f(u,z) =
(x(u),y(u),z)$ where $(u,d)\in\SS\times\RR$. The curve $\gamma(u)
= (x(u),y(u))$ is called the profile curve. The ribbon on the left is embedded, while the ribbon on the right is immersed.
The profile curve of the ribbon on the right is the self-similar figure eight curve we introduce in Section 6.}
\end{figure}

Let us briefly remark on the above procedure with respect to the particular
case of a M\"obius band.
Extending the width of the M\"obius band to infinity must result in
self-intersections, and so, the profile curve does not project uniquely to a
plane curve in $\RR^2$.
This is an essential characteristic of the M\"obius band.
One way in which we could overcome this difficulty is instead to immerse the
M\"obous band in $\R^4$, then extend the width to infinity, and project to
$\R^3$.
This would produce a profile curve in $\R^3$, the analysis of the dynamics of
which are much more complicated than for the planar curve case.

In this paper we concentrate on the dynamics of solutions to the curve
diffusion flow which evolve either by scaling, translation or rotation only.
These solutions are known as \emph{solitons}.
Solitons are of interest both from a theoretical and practical viewpoint
due to their relationship with singularities: given a family of evolving
curves, if the family exists for at most a finite amount of time, then it must
develop a singularity at the final time.
Upon rescaling the flow, effectively `zooming in on the singularity', we expect
to construct a family of curves which approximate one of the solitons for the
flow.

The existence of finite-time singularities is typically a difficult question,
especially in the case of fourth and higher-order curvature flows.
For the curve diffusion flow, this was settled in Polden \cite{Polden96}.
The analysis there was later extended by Escher and Ito \cite{EI} to more
complex initial configurations.
In particular, they prove the following theorem.

\begin{theorem}[Escher-Ito-Polden]
\label{TT}
Let $\gamma_0$ be a plane curve with winding number zero.
Let $T$ denote the maximal time of existence for the curve diffusion flow $\gamma_t := \gamma(\cdot,t)$
with initial data $\gamma_0 = \gamma(\cdot,0)$.
Then
\begin{enumerate}
\item[(i)] The estimate $T^*(\gamma_0) := \frac{L^4(\gamma_0)}{64\pi^4} \ge T$
holds, where $L(\gamma_0)$ denotes the length of $\gamma_0$; and
\item[(ii)]
If $T=T^*(\gamma_0)$ then the solution shrinks to a point as $t$ approaches $T$.
\end{enumerate}
\end{theorem}
This theorem raises two interesting questions:
\begin{enumerate}
\item[(1)] Is there a solution with winding number zero whose maximal existence
time satisfies $T=T^*(\gamma)$?
\item[(2)] Is there a solution with winding number zero that shrinks to a
point? (That is, no other kind of singularity occurs.)
\end{enumerate}
The first question investigates the sharpness of the theorem.
%%  <REVISION>
Its answer is no, since in the proof of Theorem \ref{TT}, Fenchel's inequality
is used, and for plane curves with winding number zero, it is strict.
Furthermore, the upper bound from Theorem 1 may be improved by a factor of
$\frac{\pi^2}{12}$ by using the following argument.\footnote{We would like to
thank one of the anonymous referees for pointing out this improved bound.}

Let $\gamma_t$ be a family of curves evolving by curve diffusion flow.
Differentiating the length yields
\begin{equation}
\label{Levo}
L'(\gamma_t) = - \int_\gamma |\partial_s\kappa|^2d\ell\,,
\end{equation}
where $\kappa(s)$ is the curvature of $\gamma_t$ and $\gamma(t)$ is parametrised by arc length.
We use $\partial_s$ and $d\ell$ to denote the arc length derivative and element respectively.
(We refer the reader to Section 2 for more detail on the notation.)
For curves with winding number zero, the curvature $\kappa$ has zero average,
and furthermore $\Big(\int_\gamma |\kappa|\,d\ell\Big)^2 \le
\frac{\pi^2}{12}L(\gamma_t)\int_\gamma |\kappa|^2d\ell$ (this holds for any
function with zero average), so, using the standard sharp Poincar\'e inequality
we compute
\[
\bigg(\int_\gamma |\kappa|\,d\ell\bigg)^2 \le \frac{\pi^2}{12}L\int_\gamma |\kappa|^2d\ell
 \le \frac{\pi^2}{12}L\cdot\frac{L^2}{4\pi^2}\int_\gamma |\partial_s\kappa|^2d\ell
 =  \frac{L^3}{48}\int_\gamma |\partial_s\kappa|^2d\ell
\,.
\]
In the above we have omitted the argument $\gamma_t$ from $L$.
Since $\vn{\kappa}_1 := \int_\gamma |\kappa|\,d\ell \ge 2\pi$ (by a simple variant of Fenchel's theorem), we
combine this with \eqref{Levo} above to conclude
\begin{align*}
(L^4)' &= 4L^3L' = -4\vn{\partial_s\kappa}_2^2L^3
 \le -4(48\vn{\kappa}_1^2)
 \le -3\cdot4^4\pi^2
\\
\Longrightarrow
L^4(\gamma_t) &\le -3\cdot4^4\pi^2t + L^4(\gamma_0)\,,
\end{align*}
and consequently we obtain the estimate $T \le \tilde{T}(\gamma_0) := L^4(\gamma_0)/768\pi^2$, which is
better than the estimate given by Theorem 1.

In the above proof, there are three key estimates used: Fenchel's inequality,
the Poincar\'e inequality, and a sharp H\"older-type inequality for functions
with zero average.
It is not the case that the same curve satisfies equality for all three of
these and therefore it appears clear that the improved bound $T \le
\tilde{T}(\gamma_0)$ is again not optimal.
The broader question therefore on characterising the optimal upper bound on the
maximal time of existence in terms of initial length for curves with winding
number zero remains open:
\begin{enumerate}
\item[(1')] Is there a $\hat{T}$ depending only on initial length $L(\gamma_0)$
such that:
 \begin{enumerate}
 \item All curve diffusion flows with winding number zero have maximal time bounded by $\hat{T}$; and
 \item Does there exist at least one curve diffusion flow with maximal time equal
to $\hat{T}$, or, does there exist a sequence of solutions $\gamma_i$ with
initial data $\gamma_i^{(0)}$ and winding number zero whose maximal existence
time $T_i$ satisfies
$|\hat{T}(\gamma_i^{(0)})-T_i| \rightarrow 0$?
 \end{enumerate}
\end{enumerate}
In this paper we shall not answer this question.
We do however provide an answer to the second question (2) above.
The second question is interesting from a classification of finite-time
singularities perspective.
%% </REVISION>
The key issue there is to demonstrate that solutions may become singular by
simply becoming small -- a rescaling of such solutions will exist for
all time and be asymptotic to a smooth model curve.
A finite-time singularitiy of this type is said to be of Type I, in analogy
with existing work on the mean curvature flow \cite{Hu}.
Other singularities are called Type II.
Type I singularities are expected to be simpler in structure and satisfy
natural stability criteria when compared with Type II singularities.

We answer the second question in the affirmative by way of an explicit example:
a self-similarly shrinking immersed figure eight.
This figure eight is in fact the lemniscate of Bernoulli.\footnote{The authors would like to thank Hojoo Lee for pointing this out.}
Since one requires zero area to have the possibility that the solution shrinks
to a point (see Lemma \ref{T:Marea}), and the maximal existence time is finite
(see Theorem 1), a figure eight with one leaf smaller than the other will
certainly become singular with non-zero area, ruling out the possibility that
it is asymptotic to a point.

Let $\gamma$ be the shrinking figure eight described in Section 6.
Upon scaling $\gamma$ to satisfy $L(\gamma) = 1$, we compute the
extinction time of the flow to be
\[
T = \frac{1}{3\cdot2^{11}\cdot \SK(-1)^4}\,,
\]
where $\SK(m)$ is the complete elliptic integral of the first kind with
parameter $m = k^2$.
%% <REVISION>
%This proves that the estimate $T^* \ge T$ of Theorem \ref{TT} is not sharp.
%In this particular case it overshoots $T$ by a factor of roughly three:
%% </REVISION}
In this particular case the maximal time estimate $T^*$ overshoots $T$ by a factor
of roughly 2.9, and $\tilde{T}$ overshoots $T$ by a factor of roughly 2.4:
\[
\frac{T^*}{T} = 2.9115257845\ldots\,,\qquad
\frac{\tilde{T}}{T} = 2.3946339747\ldots\,.
\]
In a sense, the self-similar figure eight represents a ``best-case'' scenario.
The generic case is that an arbitrary smooth curve with winding number zero
will be driven to a curvature singularity and not shrink to a point.
Furthermore, such a singularity involves blowup of the curvature and its
derivatives, so in particular the speed of the flow is greater in this
circumstance.
These heuristics lead us to make the following conjecture.

\begin{conj}
The maximal existence time of a curve diffusion flow emanating from a curve
with winding number zero is bounded by the existence time of a self-similar
figure eight with the same initial length.
\end{conj}
Of course, verifying the conjecture has uniqueness of the self-similar
figure eight as a corollary.
%% <REVISION>
%Although this does not answer `no' for question (1) above, it does provide some
%evidence towards a negative answer.
%The proof of Theorem \ref{TT} relies on showing that if this existence time were to
%be reached, the length of the evolving curves would be zero: They would have
%shrunk to a point.
%It is therefore natural that the best hope for achieving this bound lies in
%finding solutions which do exactly this.
%Nevertheless the question remains interesting and open.
%% </REVISION}

Another popular fourth-order surface evolution equation is the Willmore flow,
whose one-dimensional counterpart, the elastic flow, also models cylindrical
elastic ribbons.  The key difference between the two flows is that the elastic
flow minimises the integral of the square of the curvature of the curve, while
the curve diffusion flow minimises length, keeps enclosed volume fixed, and
keeps the quantity
\[
\SQ(t) = \int_0^t \bigg( \int_\gamma |\partial_s\kappa|^2\,d\ell\bigg)\,dt
\]
uniformly bounded, with a bound dependent only on the initial data.
A proof of these statements can be found in \cite{Wheeler}.  We remark on the
similarities and differences between this flow and the curve diffusion flow
throughout this article.

The structure and specific contributions of this article are as follows.  
We begin Section 2 by setting up notation, defining key geometric quantities
and giving the fourth-order nonlinear parabolic system of partial differential
equations for the curve diffusion flow.  We then briefly review relevant
previous work and the question of local existence before moving on to a
classification of stationary solutions -- those that remain constant under the
flow. We provide a five-parameter family of parametrisations which completely
classifies these solutions.  In Section 3 we record some elementary properties
of the flow that are relevant to our investigations.  In Section 4 we focus on
curves that evolve self-similarly, deriving the corresponding fourth-order
ordinary differential equation satisfied by such curves. We prove
several properties of self-similar solutions.  Section 5 is devoted to curves
that evolve by translation only. We prove there that the only closed
translators are circles, and that the only open translators (open immersions of
$\RR$) satisfying either of two conditions at infinity are straight lines.  In
Section 6 we provide an explicit parametrisation of a figure eight immersed
curve and verify that this curve does indeed evolve self-similarly under the
curve diffusion flow, contracting to a point in finite time.  In Section 7 we
consider the case of rotating solutions and prove some non-existence results
analogous to those for translators in Section 5.

\section*{Acknowledgements}

The research of the second author was supported by a summer vacation
scholarship of the Australian Mathematical Sciences Institute.  The research of
the third, fourth and fifth authors was supported by Discovery Project grant
DP120100097 of the Australian Research Council.

The authors would like to thank the two anonymous referees for their careful
reading and comments that have led to improvements in the article.
The authors would also like to thank Hojoo Lee for enlightening discussions related
to this work.

\section{The curve diffusion flow equation}

We will parametrise closed evolving curves over the unit circle $\S$, and open evolving curves over the real line $\R$.
Let $\gamma$ denote the evolving curve, with components
$$\gamma\left( u\right) = \left( x\left( u \right), y\left( u \right) \right) \mbox{.}$$
The tangent direction to the curve is given by
$$\gamma_u\left( u\right) = \left( x_u\left( u \right), y_u\left( u \right) \right) \mbox{,}$$
and therefore a choice of unit normal is
$$\nu\left( u\right) = \frac{1}{\left| \gamma_u \left( u \right)\right|} \left( y_u\left( u\right), -x_u\left( u\right) \right) \mbox{,}$$
where $\left| \cdot \right|$ denotes the length of a vector.  We will always
assume our curves are smooth and regular, that is, $\left| \gamma_u \right| >0$
for all $u$.
Recall that, in a general parametrisation, the arc length of a curve $\gamma$,
beginning at point
$\gamma\left( u_0 \right)$, is given by
\begin{equation} \label{E:length}
  s\left( u\right) = \int_{u_0}^{u} \left| \gamma_u\left( \tilde u \right) \right| d\tilde u \mbox{.}
\end{equation}
The arc length induces an arc length element $d\ell = |\gamma_u|\,du$.
The curvature vector and curvature of a plane curve is given by
%$$\kappa\left( u \right) = \frac{\left| \gamma_{uu} \times \gamma_{u} \right|}{\left| \gamma_u \right|^3} \mbox{,}$$
%where $\times$ denotes the cross product of the vectors $\gamma_u$ and $\gamma_{uu}$, where the vectors are thought of as having zero third component.
\[
\vec{\kappa} = \gamma_{ss}\quad\text{and}\quad\kappa = \ip{\nu}{\vec{\kappa}}\,,
\]
where we use $\ip{\cdot}{\cdot}$ to denote the standard inner product between vectors in the plane.
For closed embedded curves, we will always choose an inward-pointing unit normal, so that the curvature of a circle is positive.
Otherwise, we select either normal -- the curve diffusion flow is invariant under changes in orientation.

We say $\gamma_t = \gamma(\cdot,t)$ evolves under the curve diffusion flow if it satisfies the system of fourth-order partial differential equations
\begin{equation}
  \frac{\partial \gamma}{\partial t} \left( u, t\right) = - \kappa_{ss}\left( u, t\right) \nu\left( u, t\right) \quad\mbox{for}\quad \left( u, t\right) \in \mathbb{S} \times \left( 0, T\right) \label{E:CDF}
\end{equation}
with initial condition
$$\gamma\left( u, 0 \right) = \gamma_0\left( u\right) \mbox{,}$$
for some given initial curve immersion $\gamma$.  Note that an immersed curve
may have intersections, however, it has a well-defined tangent vector for each
$u$.  We label the maximal existence time of a solution to \eqref{E:CDF} as
$T$.
%We will restrict ourselves further to four times differentiable regular plane
%curves so that we may seek classical solutions to the curve diffusion flow.
%Above $\kappa$ is the curvature of curve $\gamma$ and $\nu$ is a smooth
%choice of unit normal.  We will choose $\nu$ to correspond to the outer unit
%normal for an initially embedded curve.
%The flow speed $\kappa_{ss}$ is the
%second derivative of curvature with respect to arc length.

We are interested in the shape of the evolving curves, that is, the image
$\gamma(\SS,t)$ for closed solutions, or $\gamma(\R,t)$ for open solutions, as
a set in $\R^2$, and not in the orbit of a specific point.
There is a degeneracy in the flow equation \eqref{E:CDF}: reparametrisation and
tangential movement of points in the direction of $\gamma_s$ do not affect the
shape of the image $\gamma(\SS,t)$.
The degeneracy consists of every diffeomorphism $\varphi:\SS\rightarrow\SS$.
This implies that, in terms of the images $\gamma(\SS,t)$, the flow
\eqref{E:CDF} is equivalent to
\begin{equation}
  \ip{\frac{\partial \gamma}{\partial t} \left( u, t\right)}{\nu(u,t)} = - \kappa_{ss}\left( u, t\right) \quad\mbox{for}\quad \left( u, t\right) \in \mathbb{S} \times \left( 0, T\right) \label{E:CDFN}
\end{equation}
The impact of this observation on our analysis here is that when considering
solutions with a particular form, we must also take into account a possible
diffeomorphism $\varphi$ acting on $\SS$.
This typically forces us to take normal projections in order to remove the
influence of $\varphi$.

Note that
\[
\kappa_{ss} = \partial_s^2\ip{\gamma_{ss}}{\nu}
 = \frac{1}{|\gamma_u|}\Bigg(\frac{1}{|\gamma_u|}\bigg(\ip{\frac{1}{|\gamma_u|}\Big(\frac{\gamma_u}{|\gamma_u|}\Big)_u}{\nu}\bigg)_u\Bigg)_u
\]
so that the leading order term of $-\kappa_{ss} = \ip{\partial_t\gamma}{\nu}$ is precisely
\[
-\frac{1}{|\gamma_u|^4}\ip{\gamma_{uuuu}}{\nu}\,.
\]
The normal direction is `positive' for the flow \eqref{E:CDF}, and so this
shows that \eqref{E:CDF} is a quasilinear parabolic system\footnote{It may be
helpful for the reader to recall the biharmonic heat flow $(\partial_t +
\Delta^2)v = 0$ for a function $v$.}.

Current understanding of the curve diffusion flow is quite limited, apart from
short time existence.  There are several methods for obtaining short time
existence of solutions to \eqref{E:CDF} given suitable initial data; we refer
the reader to the discussion of these in \cite{Wheeler} and the references
contained therein.  Very little else is known about the behaviour of the flow.
Giga and Ito \cite{GI1} provided an example of a simple, closed, embedded plane
curve that develops a self-intersection in finite time under the flow
\eqref{E:CDF}.  They also gave an example in \cite{GI2} of a closed, strictly
convex plane curve that becomes non-convex in finite time.  Elliot and
Maier-Paape showed that curves that are initially graphical may evolve under
\eqref{E:CDF} to become non-graphical in finite time \cite{EMP}.  On the other
hand, the fourth author recently showed \cite{Wheeler} that under the curve
diffusion flow, solution curves that are initially close to a circle, in the
sense of the normalised integral of the square of the oscillation of curvature,
exist for all time and converge exponentially fast to a (possibly different)
circle.  Further, those admissible curves which are initially embedded remain
embedded.   A fundamental obstacle in the analysis of such flows and of
\eqref{E:CDF} is that, as a system of higher order partial differential
equations, arguments based upon maximum principles are in general not
available.

%Polden \cite{Polden96} and Escher and Ito \cite{EI} have shown that smooth immersed initial curves can evolve under \eqref{E:CDF} to develop curvature singularities in finite time.  In analogy with the case for second order flows (see, for example, \cite{Hu}) it is expected that an understanding of self-similar solutions of \eqref{E:CDF} may be helpful in classifying more precisely flow behaviour at singularities.  In the case of \eqref{E:CDF}, Escher, Mayer and Simonett have provided numerical results that suggest a self-similar figure eight curve could exist.  In this paper we provide an explicit parametrisation of a figure eight immersed curve which shrinks self-similarly under \eqref{E:CDF}.

Let us now move on to the classification of stationary curve diffusion flows.
We call the family of curves $\gamma(\cdot,t)$ a \emph{closed} solution if
$\gamma$ satisfies \eqref{E:CDFN} and $\gamma(u) = \gamma(u+P)$ for some $P >
0$.  The smallest such $P$ is the period of $\gamma$.

A solution is called \emph{stationary} if it does not change in time under
\eqref{E:CDFN}, that is, if $\kappa_{ss}(u,t) \equiv 0$.

\begin{lemma}
Suppose $\gamma:\SS\rightarrow\R^2$ is a closed stationary solution to the curve diffusion flow.
Then $\gamma$ is a round circle.
\end{lemma}
\noindent \textbf{Proof:}
Integrating by parts and using $\kappa_{ss} = 0$ we have
\[
0 = -\int_\gamma \kappa\,\kappa_{ss}\,d\ell = \int_\gamma \kappa_s^2\,d\ell
\]
and so $\kappa_s = 0$, that is, $\kappa$ is constant.
Therefore $\gamma$ is a smooth closed curve with constant curvature, and by the
classification of curves in $\RR^2$, any curve with constant curvature
in the plane is a straight line or round circle.
Since $\gamma$ is closed, $\gamma$ can not be a straight line, and so it must be a round circle.
\hspace*{\fill}$\Box$

The situation is more complicated for open solutions.
We define an \emph{open curve} to be a regular immersion
$\gamma:\RR\rightarrow\RR^2$ of $\RR$ for which there does not exist a pair of
intervals $I$, $J$ such that
\begin{equation}
\label{EQtr3}
I\cap J = \emptyset\,,\quad |I| = \infty\,,\quad |J| < \infty\,,\quad\text{and}\quad \gamma(I) = \gamma(J)\,.
\end{equation}
This last condition captures the notion of not allowing $\gamma$ to `close'
while still allowing infinitely many self-intersections, including possibly
intersecting on an arbitrarily large open set.

Open curves may concentrate around a particular point, such as a logarithmic spiral.
This concentration behaviour is encapsulated by the notion of \emph{properness} of the curve $\gamma:\RR\rightarrow\RR^2$.
The curve $\gamma$ is \emph{proper} if for any compact set $C\subset\RR^2$ the inverse image
\[
\gamma^{-1}(C) := \{ u\in\RR: \gamma(u) \in C \}
\]
is a compact set.

\begin{lemma}
\label{LMopenclass}
Suppose $\gamma:\RR\rightarrow\R^2$ is an open stationary solution to the curve diffusion flow.
Then $\gamma$ is either a straight line, or a similarity transformation of the standard Cornu spiral.
Up to translation and rotation such curves satisfy
\[
\kappa(s) = k_2s + k_1
\]
for a pair $(k_1,k_2) \in \R^2$.
\end{lemma}
\noindent \textbf{Proof:}
If $\kappa_{ss} = 0$, then there is a $k_2\in\RR$ such that
\[
\kappa_s(s) = k_2\,.
\]
The curvature $\kappa$ is therefore linear in $s$; that is, there exists a $k_1\in\RR$ such that
\[
\kappa(s) = k_2s + k_1\,.
\]
If $k_2 \ne 0$, then after reparametrising by $t = s + \frac{k_1}{k_2}$ we have
\[
\kappa(t) = k_2t\,.
\]
Furthermore, by scaling the curve $\gamma$ by a factor of $\sqrt{k_2}/\sqrt{\pi}$ we have
\[
\hat{\kappa}(t) = \frac{\sqrt{\pi}}{\sqrt{k_2}} \kappa(t) = \sqrt{k_2\pi} t = \pi\hat{t}\,.
\]
Therefore either $k_2 = 0$, in which case we have a straight line, or $k_2 \ne 0$ and we have a similarity transformation of the standard Cornu spiral.
\hspace*{\fill}$\Box$\\

\noindent \textbf{Remark: } Combining the above two lemmata, we see that up to
similarity transformations the only three stationary solutions are the circle,
straight line, and Cornu spiral.
% Note that by expressing the membership criterion in
%Lemma \ref{LMopenclass} in terms of the curvature, the degeneracy by rigid
%motions in $\RR^2$ is automatically accounted for.

\begin{corollary}
Suppose $\gamma:\RR\rightarrow\R^2$ is an open, stationary, properly immersed
solution to the curve diffusion flow.
Then $\gamma$ is a straight line.
\end{corollary}
\noindent \textbf{Proof:}
We apply Lemma \ref{LMopenclass} above to obtain that $\gamma$ has curvature satisfying
\[
\kappa(s) = k_2s + k_1\,.
\]
If $k_2 = 0$ then $\kappa$ is a constant and so $\gamma$ is a straight line, as required.

Suppose on the contrary that $k_2 \ne 0$.
We will show that this leads to a contradiction to the properness assumption.
%Note that we may assume $k_2 > 0$ since if not then we may traverse the curve in the opposite direction (reparametrising $s\mapsto -s$).
Consider the restriction $\hat{\gamma}$ of $\gamma$ to the interval
$(0,\infty)$ if $k_2>0$ or to the interval $(-\infty,0)$ if $k_2<0$.
The curve $\hat{\gamma}$ has strictly monotone positive curvature.
We may therefore apply the classical Tait-Kneser theorem to obtain that the
osculating circles of $\hat{\gamma}$ are pairwise disjoint and nested (see
\cite{GTT} and the references therein for a discussion of Tait's original paper
and some interesting extensions).
Recall that the osculating circle at $\gamma(s_0)$ is the best approximating circle to $\gamma(s_0)$.
It is a standard round circle with radius $1/\kappa(s_0)$ and centre
$(\gamma(s_0) + \nu(s_0)/\kappa(s_0))$.

Let $\SC(s_0)$ denote the osculating circle at $\hat{\gamma}(s_0)$ and $\SD(s_0)$
denote the disk with boundary $\SC(s_0)$.
The theorem implies that
\[
\SD(s_0) \supset \bigcup_{s\in[s_0,\infty)} \SC(s) \supset \bigcup_{s\in[s_0,\infty)}\gamma(s)
\,.
\]
Clearly $\SD(s_0)$ is compact.
The above computation shows that 
\[
\gamma^{-1}(\SD(s_0)) = [s_0,\infty)\,,
\]
which is not compact, and so $\gamma$ is not proper.
This is contrary to the hypotheses of this Corollary, and so $k_2 = 0$, as required.
\hspace*{\fill}$\Box$\\

The two-parameter family of open curves given by Lemma \ref{LMopenclass}
therefore consists of a particular family of improperly immersed curves and
straight lines.
The improperly immersed curves are Eulerian spirals, with the clothoid or Cornu spiral ($k_2 =
\pi$ and $k_1 = 0$) a standard example.
The arc length parametrisation of the clothoid is given in terms of the Fresnel $C$ and $S$ integrals:
\[
\gamma(s) = (\,C(s), S(s)\,)\,,\quad\text{where}\quad C(s) = \int_0^s
\cos\Big(\frac{\pi t^2}{2}\Big)\,dt\quad\text{and}\quad S(s) = \int_0^s
\sin\Big(\frac{\pi t^2}{2}\Big)\,dt\,.
\]
The curvature of the clothoid is $\kappa(s) = s\pi$.
%In order to produce the family $(0,k_2)$ for an arbitrary $k_2>0$, we scale the clothoid by $\frac{1}{k_2\pi}$.
%This produces a curve with curvature $\kappa(s) = k_2s$.
%More generally, let us define the modified Fresnel integrals $C(s,c_1,c_2)$ and $\hat S(s,c_1,c_2)$ by
In order to produce the whole two-parameter family, let us define the modified Fresnel integrals $C(s,c_1,c_2)$ and $S(s,c_1,c_2)$ by
\[
C(s,c_1,c_2) = \int_0^s \cos\big(c_1t + c_2t^2\big)\,dt\quad\text{and}\quad
S(s,c_1,c_2) = \int_0^s \sin\big(c_1t + c_2t^2\big)\,dt\,.
\]
The curve $\gamma:\RR\rightarrow\RR^2$ with arc length parametrisation
\[
\gamma(s) = (\,C(s,c_1,c_2), S(s,c_1,c_2)\,)
\]
then has curvature $\kappa(s) = 2sc_2 + c_1$.

\begin{figure}
\includegraphics[width=3cm]{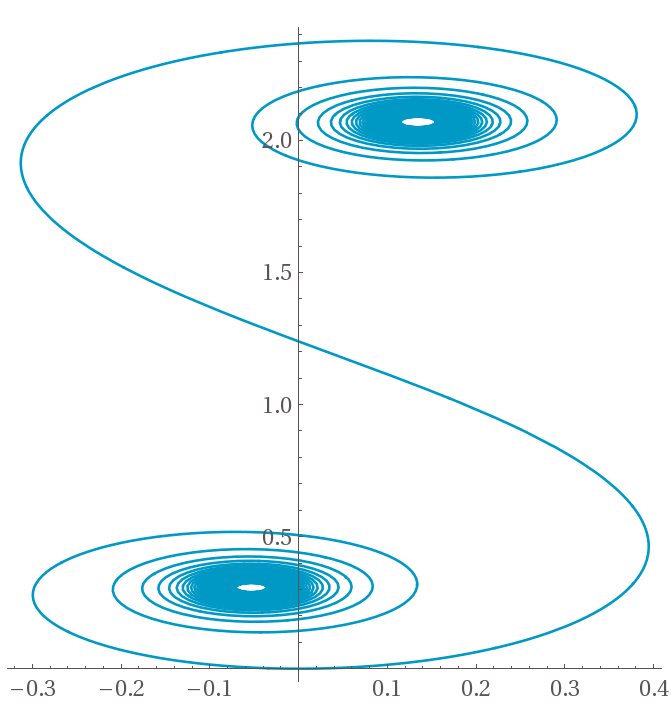}
\hspace{1cm}
\includegraphics[width=3cm]{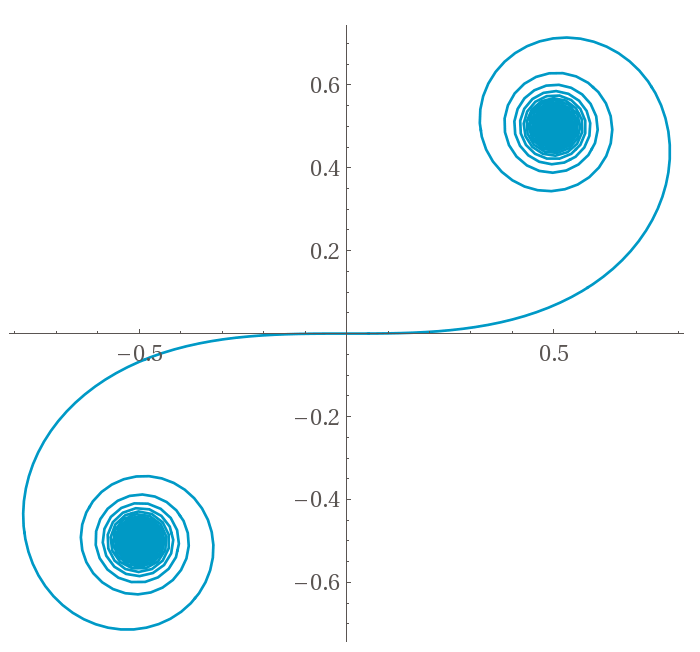}
\hspace{1cm}
\includegraphics[width=3cm]{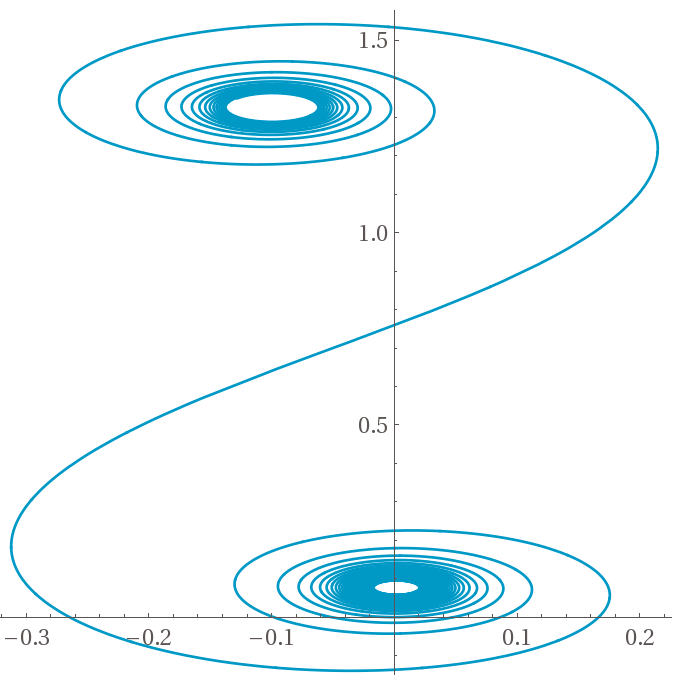}
\caption{Three members of $\SF$: the left with $\kappa(s) = -2s+3$, the middle
with $\kappa(s) = s$, and the right with $\kappa(s) = 4s+13$. We choose
$(\theta,V_1,V_2) = (0,0,0)$ for each.}
\label{FigSF}
\end{figure}

Note that this includes the case of circles ($c_2=0$) and straight lines ($c_1=c_2=0$).
The curvature $\kappa$ is invariant under rigid motions (rotations and
translations), and therefore in order to account for the entire family of
stationary solutions we must include all possible rotations and translations of
the model members $\gamma(s) = (\,C(s,c_1,c_2), S(s,c_1,c_2)\,)$.
Therefore the family $\SF$ of immersed open or closed stationary solutions to
the curve diffusion flow consists precisely of the five parameter family
\[
\SF =
\Bigg\{
  \gamma(s) = \bigg[\bigg(\begin{matrix}\cos\theta & -\sin\theta\\ \sin\theta & \cos\theta\end{matrix}\bigg)
              \bigg(\begin{matrix}C(s,c_1,c_2)\\ S(s,c_1,c_2)\end{matrix}\bigg)\bigg]^T+(V_1,V_2)\text{ for some }(c_1,c_2,V_1,V_2,\theta)\in\RR^4\times\SS
\Bigg\}\,,
\]
up to reparametrisation.
Some members of $\SF$ are depicted in Figure \ref{FigSF}.\\

\noindent \textbf{Remark} (Elastic flow)\textbf{:} We say $\gamma_t$ evolves under the elastic flow if $\gamma(\cdot,t) = \gamma_t$ satisfies
\[
\frac{\partial\gamma}{\partial t} (u,t) = \Big(-\kappa_{ss}(u,t) -\frac12\kappa^3(u,t)\Big)\nu(u,t)\,.
\]
The elastic flow is the steepest descent $L^2$-gradient flow of $\int_\gamma \kappa^2\,d\ell$.
For global regularity results and further historical references we refer to \cite{DKS}.
Stationary solutions are termed \emph{elasticae} and satisfy
\[
\kappa_{ss}(u,t) = -\frac12\kappa^3(u,t)\,.
\]
A classification of solutions to this equation is not yet known in the sense
discussed here, that is, in terms of explicit parametrisations.
Indeed, there are no known closed solutions to this equation.
Open solutions are yet to be classified.
We conjecture that there are no open properly immersed solutions apart from the
straight line.

\section{Elementary properties of the curve diffusion flow}
\newtheorem{Marea}[theorem]{Lemma}
\newtheorem{Clength}[theorem]{Lemma}

In this section we examine the behaviour of the arc length of $\gamma_t$, and the
area enclosed within $\gamma_t$, under the flow \eqref{E:CDF}.  Because $\gamma_t$
may be an immersed curve, care is required in defining the enclosed area, or
simply the change in enclosed area; we refer the reader to the discussion in
\cite{MW} for details (that article concerns evolving immersed surfaces and
discusses the definition of the signed enclosed volume, but the issues are
analogous).  When the curve $\gamma_t$ is an embedded curve, the calculations
below are the one intrinsic
dimension analogues of those in \cite{M1}, for example.
The full details for the specific case of the curve diffusion flow of immersed
curves are given in \cite{Wheeler}.

\begin{Marea} \label{T:Marea}
Under the flow \eqref{E:CDF}, the area enclosed by the curve $\gamma_t$ remains constant.
In particular, if the family of curves $\gamma_t$ contract to a point, then the initial data $\gamma_0$ must have zero signed area.
\end{Marea}

\noindent \textbf{Proof:} The signed area $A\left( \gamma_t\right)$ enclosed by $\gamma_t$ evolves under \eqref{E:CDF} according to
$$\frac{d}{dt} A\left( \gamma_t\right) = -\int_{\gamma} \kappa_{ss} \, d\ell = 0 \mbox{.}$$
That is, $A\left( t\right)$ is constant under the flow.

Since the area of a point is zero, and the enclosed area is constant under the flow, the area at any time along a family of curves that contract to a point is zero.
In particular, the area of the initial data $\gamma_0$ must be zero.
 \hspace*{\fill}$\Box$\\

\begin{Clength} \label{T:Clength}
Under the flow \eqref{E:CDF}, the length of the evolving curve $\gamma_t$ is monotone nonincreasing.
\end{Clength}

\noindent \textbf{Proof:} Under the flow \eqref{E:CDF}, the length of the evolving curve $\gamma$ evolves according to
$$\frac{d}{dt} L\left( t \right) = \int_{\gamma} \kappa_{ss} \kappa \,
d\ell = - \int_{\gamma} \kappa_s^2 \, d\ell \leq 0 \mbox{,}$$ that is, the
length is nonincreasing.
%  Moreover, the length of $\gamma$ is constant if and
%only if $\kappa_s \equiv 0$, that is, the curvature of $\gamma$ is constant.
%Since $\gamma$ is closed, this can only be the case if $\gamma$ is a
%circle.
\hspace*{\fill}$\Box$
\\

The above lemmata make the curve diffusion flow interesting from an
isoperimetric point of view.
In particular, the isoperimetric ratio is monotonically decreasing:

\begin{corollary}
The signed isoperimetric ratio
\[
I(\gamma_t) = \frac{L^2(\gamma_t)}{4\pi A(\gamma_t)}
\]
is decreasing in absolute value with value given by
\[
I(\gamma_t) = I(\gamma_0)e^{-\int_0^t\frac{2}{L(\gamma_t)}\int_\gamma \kappa_s^2\,d\ell\,d\tau}\,.
\]
\end{corollary}
\noindent \textbf{Proof:}
Differentiate $I$ and apply Lemmata \ref{T:Marea} and \ref{T:Clength} to obtain
\[
\frac{d}{dt}I(\gamma_t) = -\frac{2I(\gamma_t)}{L(\gamma_t)}\int_\gamma \kappa_s^2\,d\ell\,.
\]
Integrating yields the result.
\hspace*{\fill}$\Box$
\\

\noindent \textbf{Remark} (Elastic flow)\textbf{:} 
Although the elastic flow reduces $\int_\gamma\kappa^2\,d\ell$ with velocity
\[
\frac{d}{dt}\int_\gamma\kappa^2\,d\ell = -2\int_\gamma\Big|\kappa_{ss}+\frac12\kappa^3\Big|^2\,d\ell\,,
\]
the behaviour of the length, area, and isoperimetric ratio, are for generic
initial data unknown.

\section{Curves evolving homothetically} 
\newtheorem{SSC}[theorem]{Lemma}
\newtheorem{SSCD}[theorem]{Lemma}

A self-similar solution to a curvature flow equation is a solution whose image
maintains the same shape as it evolves, that is, it changes in time only by
scaling, translation and/or rotation.
In this section we focus on the case where $\gamma$ evolves by scaling, and in
the next section we study solutions which evolve by translation.
Section 7 studies curves evolving by rotation.

Closed, self-similar solution hypersurfaces of positive mean curvature for the
mean curvature flow are known to be spheres \cite{Hu}; the third author
generalised this result for a class of fully nonlinear curvature flows in
\cite{M4}.  In the corresponding case of the curve shortening flow,
self-similar curves were classified in \cite{AL}.
%Self-similar solutions are
%important as singularity models for various curvature flows.

\begin{definition}
We say $\gamma_t$ is a solution of \eqref{E:CDF} (or equivalently \eqref{E:CDFN})
evolving \emph{homothetically} if and only if there is a differentiable
function of time $f$ and an initial curve $\gamma_0$ such that
\[ 
\gamma\left( u, t\right) = f\left( t\right) \gamma_0\left(u\right)\,.
\]
\end{definition}

This is equivalent to
\[
\ip{\partial_t\gamma}{\nu} = f'\ip{\gamma_0}{\nu}\,.
\]
In general, homothetic solutions may be expanding, shrinking, stationary, or exhibit more complex breathing behaviour.
For the case of curve diffusion flow, Lemmata \ref{T:Marea} and \ref{T:Clength}
show that any solution must decrease its length monotonically (and be
asymptotic to a curve with constant curvature unless it exists only for a
finite amount of time). We therefore term a homothetic solution $\gamma$ a
\emph{shrinker} or \emph{shrinking solution}.
We take $f\left( 0\right) = \rho$ where $\rho\in\R$ as initial condition.
This freedom will be useful when we estimate the maximal existence time of the
shrinking figure eight in Section 6.

%, \cite{Wh}, Alexander Polden previously showed that any curve with zero winding number will either develop a singularity or shrink to a point in finite time \cite{Polden96}.  The latter fact is not obvious for a higher-than-second order parabolic partial differential equation such as \eqref{E:CDF} since the maximum principle which facilitates comparison with an enclosing curve under second order flow is no longer available.  Escher and coauthors have numerically supported the existence of a self-similar curve, but did not propose an explicit parametrisation.

The simplest example of a closed shrinking curve which is a solution to \eqref{E:CDF} is the circle.

\begin{SSC} \label{T:circ}
Let $\gamma$ be defined by $\gamma\left( u, t\right) = R\left( \cos u, \sin u\right)$ for some constant $R>0$.  Then $\gamma$ is a solution to \eqref{E:CDF}.
\end{SSC}

\noindent \textbf{Proof:}  We simply observe
$$\frac{\partial \gamma}{\partial t}\left( u, t\right) = \frac{\partial}{\partial t} R\left( \cos u, \sin u\right) = 0$$
and since the curvature $\kappa$ of $\gamma$ is identically equal to $\frac{1}{R}$, 
$$\kappa_{ss}\left( u, t\right) \equiv 0 \mbox{.}$$
Hence \eqref{E:CDF} is trivially satisfied. \hspace*{\fill}$\Box$\\

%\noindent \textbf{Remark:} If we allow smooth curves that are not closed there is a couple of other simple examples of stationary self-similar curves under \eqref{E:CDF}.  Any static straight line has zero curvature so clearly satisfies \eqref{E:CDF}, further, a static Euler spiral, or clothoid, whose curvature satisfies 
%$$\kappa\left( s\right) = \alpha s$$
%also satisfies \eqref{E:CDF}. %check, and look for quadratic in s with constant change in time?
%  \item It is interesting to consider whether \eqref{E:CDF} can have any self-similar solutions that evolve purely by translation.  Such a curve would satisfy
%  $$\gamma\left( u, t\right) = \gamma_0\left( u \right) + V t$$
%  for some constant nonzero vector $V$.  In this case
%  $$\frac{\partial \gamma}{\partial t} = V \mbox{,}$$
%  and if \eqref{E:CDF} is to hold then $\gamma$ must satisfy
%  $$-\kappa_{ss} \nu \equiv V \mbox{.}$$
%This means that the normal direction to $\gamma$ is constant and therefore the curve would have to be a straight line.  However, straight lines have zero curvature so \eqref{E:CDF} would not be satisfied.  We conclude there are no solution curves for \eqref{E:CDF} which evolve by pure translation.
%\end{enumerate}

If a curve is evolving homothetically under \eqref{E:CDF}, we may write down a
corresponding ODE by `separation of variables' that it must satisfy.

\begin{SSCD} \label{T:SSCD}
Let $\gamma: \mathbb{S} \rightarrow \mathbb{R}^2$ be a curve.  Then $\gamma$ is a self-similar solution to \eqref{E:CDF} iff there is a constant $K$ such that
\begin{equation} \label{E:SSCD}
  \kappa_{ss}\left( u\right) = K \left< \gamma\left( u\right), \nu\left( u \right) \right>
\end{equation}
for all $u\in \mathbb{S}$.
\end{SSCD}

\noindent \textbf{Proof:} Suppose the curve $\gamma$ is such that $\gamma\left(
u, t\right) = f\left( t\right) \gamma_0\left( u,t\right)$ for some
positive differentiable function $f\left( t\right)$ with $f\left( 0\right) =
\rho$, that is, at the initial time $t=0$,
$\gamma\left( u, 0\right) = \rho\,\gamma_0\left( u\right)$.
Geometric quantities associated to $\gamma$ are related to those of
$\gamma_0$ as follows.
In the calculations here we often omit the arguments of $\gamma$, $\gamma_0$ and $f$.
The tangent direction is 
\[
%\frac{\partial \gamma}{\partial u} = f\frac{\partial \gamma_0}{\partial u}
\gamma_u = f\,\partial_u\gamma_0
\]
and the unit normal is
%$$\nu\left( u, t\right) = \nu_0\left( u\right) \mbox{,}$$
\[
\nu = \nu_0\,.
\]
Above and in what follows we will include a zero subscript on geometric quantities associated with $\gamma_0$.
The curvature of $\gamma$ is related to that of $\gamma_0$ by
\begin{equation} \label{E:k}
%  \kappa\left( u, t\right) = \frac{1}{f\left( t\right) } \kappa_0 \left( u \right)\mbox{.}
\kappa = f^{-1}\kappa_0\,.
\end{equation}
In view of \eqref{E:length}, the arc length of the curves $\gamma$ and $\gamma_0$ are related via
$$s = f s_0$$
and so differentiating \eqref{E:k} twice, we find
\[
%$$\kappa_{ss}\left( u, t\right) = \frac{1}{f^3} \frac{\partial^2}{\partial s_0^2} \kappa_0\left( u\right) \mbox{.}$$
\kappa_{ss} = f^{-3}\partial^2_{s_0s_0}\kappa_0\,.
\]
Using the above we find that \eqref{E:CDF} is satisfied if and only if
\[
%$$f'\left( t\right) \gamma_0\left( u\right) = \frac{1}{f^3} \left( \kappa_0\right)_{s_0 s_0}\left( u\right) \nu_0\left( u\right) \mbox{.}$$
f'\ip{\gamma_0}{\nu_0} = -f^{-3}\partial^2_{s_0s_0}\kappa_0\ip{\nu_0}{\nu_0}\,.
\]
Following a standard separation of variables argument (note that any zeros of
$\ip{\gamma_0}{\nu_0}$ are isolated unless $\gamma$ is a straight line passing
through the origin -- this case is trivial and can be treated separately),
there must be a constant $K$ such that
\[
(f^4)' = -4\frac{\partial^2_{s_0s_0}\kappa_0}{\ip{\gamma_0}{\nu_0}} = 4K\,.
\]
Simplifying this equation finishes the proof.\hspace*{\fill}$\Box$\\

\noindent \textbf{Remarks:}
\begin{enumerate}
\item[1.] Solving the ordinary differential equation for $f$ in the above proof, we have
$$f\left( t\right) = \sqrt[4]{\rho^4 +4\, K\, t} \mbox{.}$$
Assuming also that \eqref{E:SSCD} is satisfied, we observe that
\begin{itemize}
  \item If $K=0$, then $f\left( t \right) \equiv 1$ and the curve $\gamma$ is static under \eqref{E:CDF}, that is, a stationary solution that exists for all time $t$.  
  \item If $K<0$, then the solution to \eqref{E:CDF} shrinks self-similarly to a point at time $T=-\frac{\rho^4}{4K}$.
  \item If $K>0$ then any self-similar solution would have to expand
indefinitely.  However, for closed curves this is not possible in view of Lemma
\ref{T:Clength}.  For open curves this is possible, and there could be open
expanding solutions. These could be important models for how the curve
diffusion flow smooths out isolated singularities, in the same way that the
error function generates a self-similar solution to the heat equation with the
Heaviside step function as initial data.
\end{itemize}
\item[2.]
Clearly, the shrinker condition \eqref{E:SSCD} is invariant under translation
and rotation; that is, applying these operations yields another solution to
\eqref{E:SSCD} with the same constant $K$.  Under scaling by a factor $\rho$,
one again obtains a solution to \eqref{E:SSCD}, but with a new constant $\hat K
= K\rho^4$.  One may express the ODE \eqref{E:SSCD} in terms of a pair of
equations for $(x(u), y(u))$ in the original (non-arclength) parametrisation,
and apply a standard Lie symmetry analysis to the system.
The result of this analysis is that there are no other possible transformations
that one may apply in order to obtain more solutions.
\\
\end{enumerate}

\noindent \textbf{Remark} (Elastic flow)\textbf{:} 
One naturally expects \emph{expanders} for the elastic flow, in contrast with
the shrinkers of curve diffusion flow.
The solution to the elastic flow with initial data equal to a circle expands
under the flow for all time, and so the standard round circle is the canonical
example of an expander.
Indeed, Theorems 3.2 and 3.3 of \cite{DKS} establish global existence and
convergence to an elastica under the condition that the length of the evolving
curves remain fixed.  Without this constraint, it appears impossible to obtain
convergence of the flow in any sense, since there is no way to control the
length.
There are no other known expanders for the elastic flow.

\section{Curves evolving by translation}

A family of curves $\gamma:\SS\times[0,T)\rightarrow\R^2$ evolving purely by translation satisfies
\begin{equation}
\label{E:TSCD}
\frac{\partial\gamma}{\partial t}\left( u, t\right) = \gamma_0\left( u \right) + V t + \gamma_s(u,t)\phi(u,t)
\end{equation}
for some constant vector $V$ and smooth diffeomorphism $\phi$.
%Note that $\phi$ here serves only to reparametrise the curve at each time: the image in $\RR^2$ is not affected by the action of $\phi$.
In this case
\[
\ip{\frac{\partial \gamma}{\partial t}}{\nu} = \ip{V}{\nu} \mbox{,}
\]
and if \eqref{E:CDF} is to hold then $\gamma$ must satisfy
\begin{equation}
\label{E:TSCD2}
-\kappa_{ss} \equiv \ip{V}{\nu} \mbox{.}
\end{equation}
We call the solution $\gamma$ a translator.
If $V = 0$, then the solution $\gamma$ is stationary, and we call the translator trivial.

It is common that for a given curvature flow there are no closed non-trivial translators.
Nevertheless one typically expects open immersed non-trivial translating solutions.
A family of open immersed curves $\gamma:\R\times[0,T)\rightarrow\R^2$ evolving purely by translation satisfies
\eqref{E:TSCD} and \eqref{E:TSCD2} as in the closed setting.
Such curves usually arise as blowup limits of a Type 2 singularity \cite{GageHamilton}.
Singularities of this type have engendered much interest in the literature.

In this section we investigate translators for the curve diffusion flow.
We prove by a simple argument that there are no smooth closed non-trivial translators.
Our proof implies a stronger statement about open curves, which is sharp, by the earlier example of the clothoid.

\begin{proposition}
\label{P:nonontrivialTs}
Let $\gamma:\SS\rightarrow\R^2$ be a smooth closed translator.
Then $\gamma$ is trivial; that is, $\gamma(\SS)$ is a standard round circle.
\end{proposition}
\noindent \textbf{Proof:}
Since $\gamma$ is closed, integration by parts gives
\begin{equation}
\label{EQtr1}
- \int_\gamma \kappa\kappa_{ss}\,d\ell = \int_\gamma \kappa_s^2d\ell\,.
\end{equation}
Also, noting that $V$ is a constant vector and that the curvature vector
$\vec{\kappa}$ satisfies $\kappa\nu = \vec{\kappa} = \gamma_{ss}$, we have
\begin{equation}
\label{EQtr2}
\int_\gamma \ip{V}{\vec{\kappa}}\,d\ell = \int_\gamma \partial_s\big(\ip{V}{\gamma_s}\big)\,d\ell = 0\,.
\end{equation}
Multiplying \eqref{E:TSCD2} by $\kappa$ and integrating gives
\[
- \int_\gamma \kappa\kappa_{ss}\,d\ell = \int_\gamma \ip{V}{\vec{\kappa}}\,d\ell\,,
\]
which when combined with \eqref{EQtr1} and \eqref{EQtr2} above yields
\[
\int_\gamma \kappa_s^2d\ell = 0\,.
\]
Since $\gamma$ is smooth, this implies that $\kappa_s \equiv 0$ on $\gamma$, that is, that the curvature is constant.
Therefore $\kappa_{ss} = 0$, and so $V=0$; that is, the translator is trivial.
By Lemma 2, it is a round circle.
\hspace*{\fill}$\Box$
\\

Now let us investigate open curves.

\begin{proposition}
\label{P:nonontrivialTs2}
Suppose $\gamma:\RR\rightarrow\R^2$ is a smooth open translator satisfying
\begin{equation}
\label{EQgcond}
\lim_{i\rightarrow\infty} \Big[\big(\kappa(p_i)\kappa_s(p_i) + \ip{V}{\gamma_s(p_i)}\big) - \big(\kappa(q_i)\kappa_s(q_i) + \ip{V}{\gamma_s(q_i)}\Big)\Big]= 0
\end{equation}
for arbitrary sequences $\{p_i\}$ and $\{q_i\}$ asymptotic to $-\infty$ and $+\infty$ respectively.
Then $\gamma(\RR)$ is a straight line.
\end{proposition}
\noindent \textbf{Proof:}
Let us follow the idea of Proposition \ref{P:nonontrivialTs}.
We integrate \eqref{E:TSCD2} against $\kappa$ to find
\begin{align*}
\int_{p_i}^{q_i} \kappa_s^2d\ell
  &= - \int_{p_i}^{q_i} \kappa\kappa_{ss}\,d\ell    + \kappa\kappa_s\big|_{p_i}^{q_i}
\\
  &=   \int_{p_i}^{q_i} \ip{V}{\vec{\kappa}}\,d\ell + \kappa\kappa_s\big|_{p_i}^{q_i}
\\
  &=   \int_{p_i}^{q_i} \partial_s\big(\ip{V}{\gamma_s}\big)\,d\ell + \kappa\kappa_s\big|_{p_i}^{q_i}
\\
  &=   \Big[  \ip{V}{\gamma_s} + \kappa\kappa_s\Big]_{p_i}^{q_i}\,.
\end{align*}
Taking $i\rightarrow\infty$ yields the result.
Note that circles are excluded by the open hypothesis.
\hspace*{\fill}$\Box$
\\

%\noindent \textbf{Remark:}
%It is straightforward to show that the results of Propositions \ref{P:nonontrivialTs} and \ref{P:nonontrivialTs2} extend to curves in the regularity classes $W^{3,2}(\SS)$ and $W^{3,2}_{loc}(\RR)$ respectively.

\noindent \textbf{Remark:}
The result of Proposition \ref{P:nonontrivialTs2} is false without condition \eqref{EQgcond}.
To see this, consider the clothoid: a smooth open curve with $\kappa(s) = s$ where $s$ is the arc length parameter.
Clearly this curve satisfies $\kappa_{ss} = 0$ and so is a stationary translator with $V=0$ in \eqref{E:TSCD2}.
The clothoid however does not satisfy condition \eqref{EQgcond}, since $\kappa(s)\kappa_s(s) = s$ is an odd function, and so the limit does not exist.\\

Although the above remark shows in some sense the sharpness of Proposition
\ref{P:nonontrivialTs2}, it leaves open the possibility of replacing the decay
condition \eqref{EQgcond} with something different.

The following proposition explores this and converts \eqref{EQgcond} from a
pointwise decay condition on the curvature to a natural growth condition on $\kappa$ in $L^2$ and $V$ in $L^1$.

\begin{proposition}
\label{P:nonontrivialTs3}
Suppose $\gamma:\RR\rightarrow\R^2$ is a smooth open proper translator satisfying
\begin{equation}
\label{EQgcond2}
\lim_{\rho\rightarrow\infty}\bigg( \frac1{\rho^2}\int_{\gamma^{-1}(B_\rho(0))} \kappa^2\,d\ell
   + \frac{1}{\rho}\int_{\gamma^{-1}(B_{2\rho}(0))} |\ip{V}{\gamma_s}|\,d\ell\,.
\bigg) = 0\,,
\end{equation}
where $B_\rho(0)$ denotes the closed ball of radius $\rho$ centred at the origin.
Then $\gamma(\RR)$ is a straight line.
\end{proposition}
\noindent \textbf{Proof:}
We use again a similar idea, however this time we localise the estimate in the
plane using a cutoff function $\eta_\rho = \tilde{\eta}_\rho\circ\gamma$ where
$\tilde{\eta}_\rho:\R^2\rightarrow[0,1]$ satisfies
\begin{align*}
\tilde{\eta}_\rho(x) &= 1\quad\text{on}\quad B_{\frac\rho2}(0)\,,
\\
\tilde{\eta}_\rho(x) &= 0\quad\text{on}\quad B_{\rho}(0)\,,
\\
\eta\text{ is of class $C^\infty$,}&\text{ supp $\eta$ is compact, and}
\\
\text{ there exist constants }&c_k\text{ such that }\partial_s^k\eta \le \frac{c_k}{\rho^k}\,.
\end{align*}
The existence of such a function is straightforward.
A constructive proof can be found in \cite{WillmoreRiemannianGeometry}.
Let us denote by $\gamma^{-1}(S)$ the inverse image of a set $S\subset\R^2$.
First, note that since $\gamma$ is proper, the integral
\[
\int_{\gamma^{-1}(B_\rho(0))} \kappa_s^2\,\eta^4\, d\ell
\]
is well-defined, irrespective of $\rho$.
We compute
\begin{align*}
\int_{\gamma^{-1}(B_\rho(0))} \kappa_s^2\,\eta^4\, d\ell
  &= - \int_{\gamma^{-1}(B_\rho(0))} \kappa\kappa_{ss}\,\eta^4\,d\ell
     - 4\int_{\gamma^{-1}(B_\rho(0))} \kappa\kappa_s\,\eta_s\,\eta^3\,d\ell
\\
  &=   \int_{\gamma^{-1}(B_\rho(0))} \ip{V}{\vec{\kappa}}\,\eta^4\,d\ell
     - 4\int_{\gamma^{-1}(B_\rho(0))} \kappa\kappa_s\,\eta_s\,\eta^3\,d\ell
\\
  &=   \int_{\gamma^{-1}(B_\rho(0))} \partial_s\big(\ip{V}{\gamma_s}\big)\,\eta^4\,d\ell
     - 4\int_{\gamma^{-1}(B_\rho(0))} \kappa\kappa_s\,\eta_s\,\eta^3\,d\ell
\\
  &= - 4\int_{\gamma^{-1}(B_\rho(0))} \ip{V}{\gamma_s}\,\eta_s\,\eta^3\,d\ell
     - 4\int_{\gamma^{-1}(B_\rho(0))} \kappa\kappa_s\,\eta_s\,\eta^3\,d\ell
\\
  &\le
       \varepsilon\int_{\gamma^{-1}(B_\rho(0))} \kappa_s^2\,\eta^4\, d\ell
     + 4\int_{\gamma^{-1}(B_\rho(0))} \big|\ip{V}{\gamma_s}\big|\,|\eta_s|\,\eta^3\,d\ell
\\&\qquad
     + \frac{4}{\varepsilon}\int_{\gamma^{-1}(B_\rho(0))} \kappa^2\,\eta_s^2\,\eta^2\,d\ell\,,
\end{align*}
for any $\varepsilon>0$.
Choosing $\varepsilon = \frac14$ and absorbing the first term on the left, we find
\[
\int_{\gamma^{-1}(B_\rho(0))} \kappa_s^2\,\eta^4\, d\ell
  \le
     \frac{c}{\rho^2}\int_{\gamma^{-1}(B_{2\rho}(0))} \kappa^2\,d\ell
   + \frac{c}{\rho}\int_{\gamma^{-1}(B_{2\rho}(0))} |\ip{V}{\gamma_s}|\,d\ell\,.
\]
Taking $\rho\rightarrow\infty$ in this estimate yields $\kappa_s = 0$ and hence the result.
\hspace*{\fill}$\Box$
\\

{\bf Remark:} It is an open problem to determine if either the properness
condition or growth condition \eqref{EQgcond2} are necessary in Proposition
\ref{P:nonontrivialTs3}.
One advantage of the integral growth condition over the pointwise decay
condition is that the rigidity statement holds for a much wider class of weak
solutions, which may not even have continuous tangent vector.
We will not explore notions of weak solution here.

\section{The shrinking figure eight}
\newtheorem{figure8}[theorem]{Proposition}

\begin{figure8} \label{T:figure8}
Let $\gamma: \mathbb{S} \rightarrow \mathbb{R}^2$ be defined by
\begin{equation} \label{E:figure8}
  \gamma\left( u \right) = \frac{1}{\left( 1 + \sin^2 u \right)} \left( \cos u, \frac{1}{2} \sin 2u \right) \mbox{.}
\end{equation}
Then $\gamma$ evolves self-similarly under \eqref{E:CDF}.
\end{figure8}

{\bf Remark:} The curve defined by \eqref{E:figure8} is the lemniscate of Bernoulli.

\noindent \textbf{Proof:} The tangent direction to $\gamma$ is given by
$$\gamma_u \left( u \right) = \frac{1}{\left( 2 - \cos^2 u\right)^2} \left( - \left( 2 + \cos^2 u\right) \sin u, -2 + 3 \cos^2 u \right)$$
and the unit normal is
$$\nu\left( u\right) = \left( \frac{3 \sin^2 u - 1}{\left( 1 + \sin^2 u\right)^{\frac{3}{2}}}, \frac{-\sin u \left( 2 + \cos^2 u \right)}{\left( 1 + \sin^2 u\right)^{\frac{3}{2}}} \right) \mbox{.}$$
Therefore, we have
\begin{equation} \label{E:inner}
  \left< \gamma\left( u\right), \nu\left( u\right) \right> = \frac{\cos u\left( \sin^2 u - 1\right)}{\left( 1 + \sin^2 u\right)^{\frac{3}{2}}} \mbox{.}
  \end{equation}

By further straightforward calculations, the curvature of $\gamma$ is given by
$$\kappa\left( u\right) = \frac{3 \cos\left( u\right)}{\left( 1 + \sin^2 u \right)^{\frac{1}{2}}} \mbox{.}$$
We calculate using \eqref{E:length} with the chain rule that the derivatives or curvature with respect to arc length are 
$$\kappa_s\left( u\right) = \frac{-6 \sin u}{1 + \sin^2 u} \quad\mbox{and}\quad
\kappa_{ss}\left( u \right) = \frac{-6 \cos^3 u}{\left( 1 + \sin^2
u\right)^{\frac{3}{2}}} \mbox{.}$$

Comparing with \eqref{E:inner}, we see that
$$6\left< \gamma\left( u \right) , \nu\left( u \right) \right> = \kappa_{ss}\left( u \right) \mbox{,}$$
so $\gamma$ given by \eqref{E:figure8} satisfies \eqref{E:SSCD} with $K= - 6$.\hspace*{\fill}$\Box$\\

\noindent \textbf{Remarks:}
\begin{enumerate}
\item[1.] Since $K<0$ above, in view of the Remark at the end
of Section 4, the figure eight given by \eqref{E:figure8} shrinks to a point at
time $T=\frac{1}{24}$.  This is in fact consistent with Lemma \ref{T:Marea},
since the symmetry of the figure eight implies that the signed area is equal to
zero and remains constant under \eqref{E:CDF}.  More generally, in view of
Lemma \ref{T:Marea}, any curve that shrinks self-similarly under \eqref{E:CDF}
must have zero signed area.
\item[2.] The length of the parametrisation of the figure eight above is
\[
L(\gamma) = 4\SK(-1)\,.
\]
Therefore, choosing $\rho = \frac1{4\SK(-1)}$ allows us to obtain the maximal
existence time of the figure eight with unit length as claimed in the
introduction.
\end{enumerate}

\begin{figure}
\includegraphics[width=3cm]{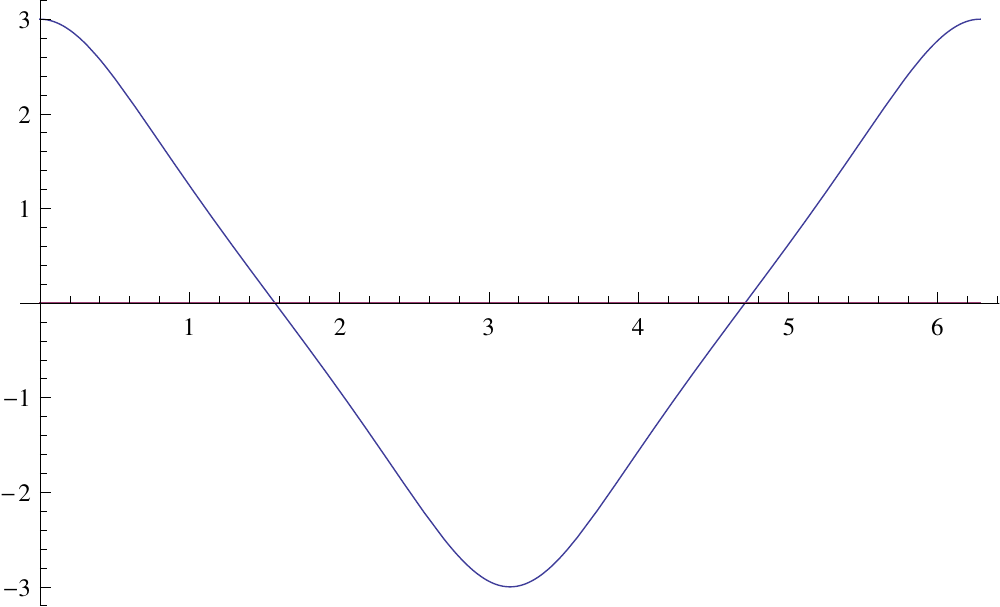}
\hspace{1cm}
\includegraphics[width=5cm]{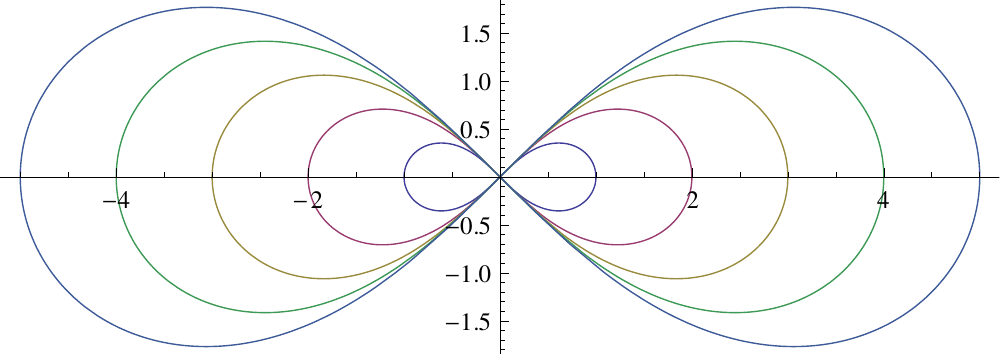}
\hspace{1cm}
\includegraphics[width=3cm]{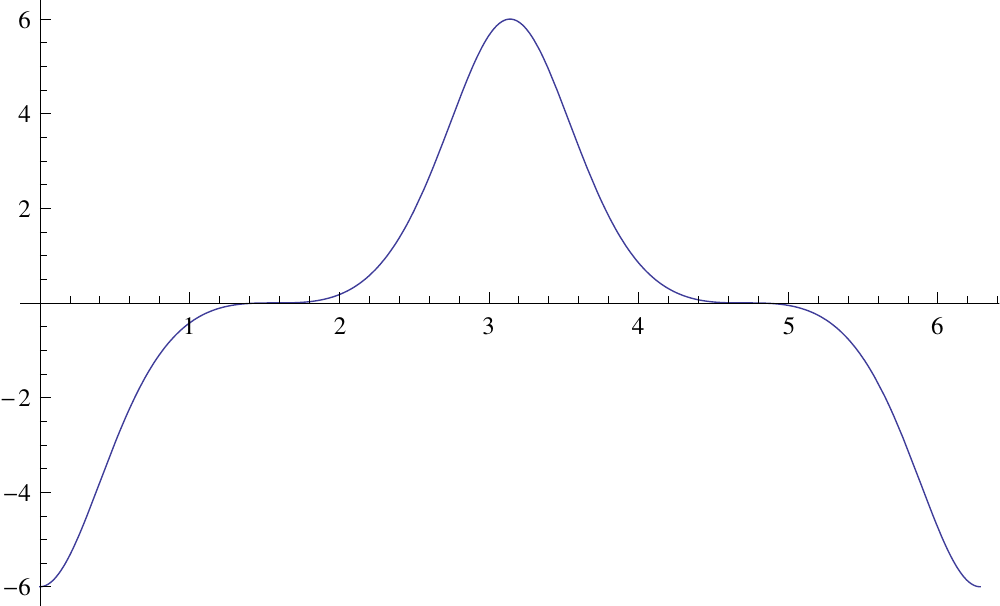}
\caption{The shrinking figure eight (middle), along with plots of its curvature
(left) and second arc length derivative of curvature (right).
The horizontal axis represents the $u$ variable and the vertical axis
represents $\kappa(u)$ and $\kappa_{ss}(u)$ in the left and right figures
respectively.}
\end{figure}

\section{Curves evolving by rotation}

A family of curves $\gamma:\SS\times[0,T)\rightarrow\R^2$ evolving purely by rotation satisfies
\begin{equation}
\label{E:RSCD}
-\kappa_{ss} \equiv 2S(t)\ip{\gamma_s}{\gamma} \mbox{.}
\end{equation}
In particular, if \eqref{E:RSCD} is satisfied for a curve $\gamma_0$, then the
curve diffusion flow evolving from that curve evolves purely by translation.
We call the solution $\gamma_t$ a rotator.

Rotators have been conjectured to exist for the curve shortening flow for some
time, with a specific conjecture in Altschuler \cite{A91}.
A rigorous classification of these has only recently appeared in the literature
however \cite{H12}.
In this section we prove some classification results for rotators.
Our method is to multiply \eqref{E:RSCD} by the curvature and integrate.
%One might expect that integrating \eqref{E:RSCD} without multiplying by
%anything is a better strategy, since both sides are exact (i.e. the equation
%reads $(\kappa_{s} + S|\gamma|^2)_s = 0$).
%Proofs beginning with this idea seem to be longer and more complicated than
%multiplying with $\kappa$, however, and so we have used this approach below.

\begin{proposition}
\label{PRnonontrivialRs}
Let $\gamma:\SS\rightarrow\R^2$ be a smooth closed rotator.
Then $\gamma$ is a standard round circle.
\end{proposition}
\noindent \textbf{Proof:}
The unit normal vector satisfies
\[
0 = \frac12\partial_s|\nu|^2 = \ip{\nu}{\nu_s}\,.
\]
%and since $\{\gamma_s,\nu\}$ is an orthonormal basis of $\R^2$, we have
%\begin{equation}
%\label{EQnus}
%\nu_s = \ip{\nu_s}{\gamma_s}\gamma_s = -\ip{\nu}{\gamma_{ss}}\gamma_s = -\kappa\gamma_s\,,
%\end{equation}
%where we used $0 = \partial_s\ip{\gamma_s}{\nu} = \ip{\gamma_{ss}}{\nu} + \ip{\gamma_s}{\nu_s}$ for the second equality.
%Noting that $2S\ip{\gamma_s}{\gamma} = S\partial_s|\gamma|^2$ and using \eqref{EQnus} above,
%Using \eqref{EQnus} above,
Using the Frenet-Serret frame equations,
\begin{equation*}
\int_\gamma \kappa_s^2\,d\ell
 = - \int_\gamma \kappa\,\kappa_{ss}\,d\ell
 = -2S \int_\gamma \ip{k\gamma_s}{\gamma}\,d\ell
 =  2S \int_\gamma \ip{\nu_s}{\gamma}\,d\ell
 = -2S \int_\gamma \ip{\nu}{\gamma_s}\,d\ell
 = 0\,,
\end{equation*}
irrespective of $S$.
The curvature of $\gamma$ is therefore constant and its image must be a round circle.
\hspace*{\fill}$\Box$
\\

Let us now turn to open curves.

\begin{proposition}
\label{PRnonontrivialRs2}
Suppose $\gamma:\RR\rightarrow\R^2$ is a smooth open rotator satisfying
\begin{equation}
\label{EQgcondr}
\lim_{i\rightarrow\infty}
  \Big[\big(\kappa(p_i)\kappa_s(p_i) - 2S\ip{\nu(p_i)}{\gamma(p_i)}\big) - \big(\kappa(q_i)\kappa_s(q_i) - 2S\ip{\nu(q_i)}{\gamma(q_i)}\Big)\Big]
= 0
\end{equation}
for arbitrary sequences $\{p_i\}$ and $\{q_i\}$ asymptotic to $-\infty$ and $+\infty$ respectively.
Then $\gamma(\RR)$ is a straight line.
\end{proposition}
\noindent \textbf{Proof:}
Integrate \eqref{E:RSCD} against $\kappa$ to find
\begin{align*}
\int_{p_i}^{q_i} \kappa_s^2d\ell
  &= - \int_{p_i}^{q_i} \kappa\kappa_{ss}\,d\ell    + \kappa\kappa_s\big|_{p_i}^{q_i}
\\
  &= - 2S\int_{p_i}^{q_i} \ip{\nu_s}{\gamma}\,d\ell + \kappa\kappa_s\big|_{p_i}^{q_i}
\\
  &=   \Big[- 2S\ip{\nu}{\gamma} + \kappa\kappa_s\Big]_{p_i}^{q_i}\,.
\end{align*}
Taking $i\rightarrow\infty$ yields the result.
\hspace*{\fill}$\Box$
\\

\noindent \textbf{Remark:}
The result of Proposition \ref{PRnonontrivialRs2} is false without condition
\eqref{EQgcondr}, for the same reason that Proposition \ref{P:nonontrivialTs2}
is false without condition \eqref{EQgcond}.
%  For the convenience of the reader
%we provide it again here.
%Consider the clothoid: a smooth open curve with $\kappa(s) = s$ where $s$ is the arc length parameter.
%Clearly this curve satisfies $\kappa_{ss} = 0$ and so is a stationary rotator
%with $S=0$ in \eqref{E:TSCD2}.
%The clothoid however does not satisfy condition \eqref{EQgcondr}, since
%$\kappa(s)\kappa_s(s) = s$ is an odd function, and so the limit does not
%exist.
\\

As with translators earlier, we may convert \eqref{EQgcondr} from a
pointwise decay condition on the curvature to a natural growth condition on
$\kappa$ in $L^2$ and $S\ip{\nu}{\gamma}$ in $L^1$.

\begin{proposition}
\label{PRnonontrivialRs3}
Suppose $\gamma:\RR\rightarrow\R^2$ is a smooth open proper rotator satisfying
\begin{equation}
\label{EQgcond2r}
\lim_{\rho\rightarrow\infty}\bigg( \frac1{\rho^2}\int_{\gamma^{-1}(B_\rho(0))} \kappa^2\,d\ell
   + \frac{|S|}{\rho}\int_{\gamma^{-1}(B_{2\rho}(0))} \big|\ip{\nu}{\gamma}\big|\,d\ell
\bigg) = 0\,,
\end{equation}
where $B_\rho(0)$ denotes the closed ball of radius $\rho$ centred at the origin.
Then $\gamma(\RR)$ is a straight line.
\end{proposition}
\noindent \textbf{Proof:}
Let $\eta_\rho = \tilde{\eta}_\rho\circ\gamma$ be the cutoff function used in the proof of Proposition \ref{P:nonontrivialTs3}.
We compute
\begin{align*}
\int_{\gamma^{-1}(B_\rho(0))} \kappa_s^2\,\eta^4\, d\ell
  &= - \int_{\gamma^{-1}(B_\rho(0))} \kappa\kappa_{ss}\,\eta^4\,d\ell
     - 4\int_{\gamma^{-1}(B_\rho(0))} \kappa\kappa_s\,\eta_s\,\eta^3\,d\ell
\\
  &= 2S\int_{\gamma^{-1}(B_\rho(0))} \ip{\nu_s}{\gamma}\,\eta^4\,d\ell
     - 4\int_{\gamma^{-1}(B_\rho(0))} \kappa\kappa_s\,\eta_s\,\eta^3\,d\ell
\\
  &= - 8S\int_{\gamma^{-1}(B_\rho(0))} \ip{\nu}{\gamma}\,\eta_s\,\eta^3\,d\ell
     - 4\int_{\gamma^{-1}(B_\rho(0))} \kappa\kappa_s\,\eta_s\,\eta^3\,d\ell
\\
  &\le
       \varepsilon\int_{\gamma^{-1}(B_\rho(0))} \kappa_s^2\,\eta^4\, d\ell
     + 8|S|\int_{\gamma^{-1}(B_\rho(0))} \big|\ip{\nu}{\gamma}\big|\,|\eta_s|\,\eta^3\,d\ell
\\&\qquad
     + \frac{4}{\varepsilon}\int_{\gamma^{-1}(B_\rho(0))} \kappa^2\,\eta_s^2\,\eta^2\,d\ell\,,
\end{align*}
for any $\varepsilon>0$.
Choosing $\varepsilon = \frac14$ and absorbing the first term on the left, we find
\[
\int_{\gamma^{-1}(B_\rho(0))} \kappa_s^2\,\eta^4\, d\ell
  \le
     \frac{c}{\rho^2}\int_{\gamma^{-1}(B_{2\rho}(0))} \kappa^2\,d\ell
   + \frac{c|S|}{\rho}\int_{\gamma^{-1}(B_{2\rho}(0))} \big|\ip{\nu}{\gamma}\big|\,d\ell\,.
\]
Taking $\rho\rightarrow\infty$ in this estimate yields $\kappa_s = 0$ and hence the result.
\hspace*{\fill}$\Box$
\\

\end{document}